\documentclass[fleqn,11pt,twoside]{article}

\usepackage{amsthm,amsthm,amssymb, color, xcolor,epsfig, graphics, subfigure}

\usepackage{amsmath, graphicx, latexsym, lscape }

\makeatletter
\newcommand{\copyrightnote}[2]{{\renewcommand{\thefootnote}{}
 \footnotetext{\small\it
\begin{flushleft}
 \copyright \ #1   #2  
\end{flushleft}}}}

\newcommand{\Name}[1]{\begin{flushleft}
                       \LARGE \bf #1
                       \end{flushleft}\vspace{-3mm}}

\newcommand{\Author}[1]{\begin{flushleft}
                       \it #1 \end{flushleft}}

\newcommand{\Address}[1]{\begin{flushleft}
                       \it #1 \end{flushleft}}

\newcommand{\Date}[1]{\begin{flushleft}
                      \small  \it #1 \end{flushleft}}

\newcommand{\evenhead}{Author \ name}
\newcommand{\oddhead}{Article \ name}

\renewcommand{\@evenhead}{
\hspace*{-3pt}\raisebox{-15pt}[\headheight][0pt]{\vbox{\hbox to \textwidth
{\thepage \hfil \evenhead}\vskip4pt \hrule}}}
\renewcommand{\@oddhead}{
\hspace*{-3pt}\raisebox{-15pt}[\headheight][0pt]{\vbox{\hbox to \textwidth
{\oddhead \hfil \thepage}\vskip4pt\hrule}}}
\renewcommand{\@evenfoot}{}
\renewcommand{\@oddfoot}{}

\long\def\@makecaption#1#2{%
  \vskip\abovecaptionskip
  \sbox\@tempboxa{\small \textbf{#1.}\ \ #2}%
  \ifdim \wd\@tempboxa >\hsize
    {\small \textbf{#1.}\ \ #2}\par
  \else
    \global \@minipagefalse
    \hb@xt@\hsize{\hfil\box\@tempboxa\hfil}%
  \fi
  \vskip\belowcaptionskip}

\newcommand{\JNMPnumberwithin}[3][\arabic]{%
  \@ifundefined{c@#2}{\@nocounterr{#2}}{%
    \@ifundefined{c@#3}{\@nocnterr{#3}}{%
      \@addtoreset{#2}{#3}%
      \@xp\xdef\csname the#2\endcsname{%
        \@xp\@nx\csname the#3\endcsname .\@nx#1{#2}}}}%
}

\newcommand{\resetfootnoterule} {
  \renewcommand\footnoterule{%
  \kern-3\p@
  \hrule\@width.4\columnwidth
  \kern2.6\p@}
}

\renewcommand{\footnoterule}{}

\makeatother

\begin{document}

%\setcounter{MaxMatrixCols}{10}

%TCIDATA{OutputFilter=LATEX.DLL}
%TCIDATA{Version=5.50.0.2960}
%TCIDATA{<META NAME="SaveForMode" CONTENT="1">}
%TCIDATA{BibliographyScheme=Manual}
%TCIDATA{Created=Monday, January 18, 2021 15:26:55}
%TCIDATA{LastRevised=Saturday, May 15, 2021 08:23:20}
%TCIDATA{<META NAME="GraphicsSave" CONTENT="32">}
%TCIDATA{<META NAME="DocumentShell" CONTENT="Standard LaTeX\Blank - Standard LaTeX Article">}
%TCIDATA{Language=American English}
%TCIDATA{CSTFile=40 LaTeX article.cst}

\newtheorem{theorem}{Theorem}
\newtheorem{acknowledgement}[theorem]{Acknowledgement}
\newtheorem{algorithm}[theorem]{Algorithm}
\newtheorem{axiom}[theorem]{Axiom}
\newtheorem{case}[theorem]{Case}
\newtheorem{claim}[theorem]{Claim}
\newtheorem{conclusion}[theorem]{Conclusion}
\newtheorem{condition}[theorem]{Condition}
\newtheorem{conjecture}[theorem]{Conjecture}
\newtheorem{corollary}[theorem]{Corollary}
\newtheorem{criterion}[theorem]{Criterion}
\newtheorem{exercise}[theorem]{Exercise}
\newtheorem{lemma}[theorem]{Lemma}
\newtheorem{notation}[theorem]{Notation}
\newtheorem{problem}[theorem]{Problem}
\newtheorem{proposition}[theorem]{Proposition}
\newtheorem{remark}[theorem]{Remark}
\newtheorem{solution}[theorem]{Solution}
\newtheorem{summary}[theorem]{Summary}

\setcounter{page}{1}

\renewcommand{\evenhead}{ {\LARGE\textcolor{blue!10!black!40!green}{{\sf \ \ \ ]ocnmp[}}}\strut\hfill F Calogero and F Payandeh}
\renewcommand{\oddhead}{ {\LARGE\textcolor{blue!10!black!40!green}{{\sf ]ocnmp[}}}\ \ \ \ \
Explicitly solvable systems of two autonomous first-order ODEs
}

%%%% Matter for the first page 
\thispagestyle{empty}
\newcommand{\FistPageHead}[3]{
\begin{flushleft}
\raisebox{8mm}[0pt][0pt]
{\footnotesize \sf
\parbox{150mm}{{Open Communications in Nonlinear Mathematical Physics}\ \  \ \ {\LARGE\textcolor{blue!10!black!40!green}{]ocnmp[}}
\qquad Vol.1 (2021) pp
#2\hfill {\sc #3}}}\vspace{-13mm}
\end{flushleft}}

\FistPageHead{1}{\pageref{firstpage}--\pageref{lastpage}}{ \ \ Article}

\strut\hfill

\copyrightnote{The author(s). Distributed under a Creative Commons Attribution 4.0 International License}
%{2021}{Authors' Names} what

\Name{Explicitly solvable systems of two autonomous first-order ordinary
differential equations with homogeneous quadratic right-hand sides }

\label{firstpage}

\Author{
Francesco Calogero$^{\,a,b,1}$ and Farrin Payandeh$^{\,c,2}$
}

\Address{$^{a}$
Physics Department, University of Rome "La Sapienza", Rome,
Italy
\\[2mm]
$^{b}$
Istituto Nazionale di Fisica Nucleare, Sezione di Roma 1\\[2mm]
$^{c}$
 Department of Physics, Payame Noor University (PNU)\\[2mm]
 \ \  PO Box
19395-3697 Tehran, Iran\\[2mm]
$^{1}$ francesco.calogero@uniroma1.it;
francesco.calogero@roma1.infn.it\\[2mm]
$^{2}$ f\_payandeh@pnu.ac.ir; farrinpayandeh@yahoo.com
}

\Date{Received Date: 8 April 2021; Accepted Date: 15 May 2021}

\renewcommand{\theequation}{\arabic{equation}}
\setcounter{equation}{0}

\begin{abstract}
\noindent
After tersely reviewing the various meanings that can be given to the
property of a system of nonlinear ODEs to be \textit{solvable}, we identify
a special case of the system of two first-order ODEs with \textit{%
homogeneous quadratic} right-hand sides which is \textit{explicitly solvable}%
. It is identified by $2$\textit{\ explicit algebraic constraints} on the $6$
\textit{a priori arbitrary} parameters that characterize this system. Simple
extensions of this model to cases with \textit{nonhomogeneous quadratic}
right-hand sides are also identified, including \textit{isochronous} cases.
\end{abstract}

%\renewcommand{\theequation}{\thesection.\arabic{equation}}

%\renewcommand{\theequation}{\arabic{section}.\arabic{equation}}

%\renewcommand{\theequation}{\thesection.\arabic{equation}}

%%%%%%%%%%%%%%%%%%%%%%

 \setcounter{equation}{0}

\section{Introduction}

In this paper we mainly focus on the following system of $2$ \textit{%
first-order} ODEs with \textit{homogeneous quadratic} right-hand sides:%
\begin{equation}
\dot{x}_{n}\left( t\right) =c_{n1}\left[ x_{1}\left( t\right) \right]
^{2}+c_{n2}x_{1}\left( t\right) x_{2}\left( t\right) +c_{n3}\left[
x_{2}\left( t\right) \right] ^{2}~,~~~n=1,2~.  \label{1}
\end{equation}

\textbf{Notation 1-1}. Above and hereafter $t$ is the \textit{independent}
variable, and superimposed dots indicate \textit{differentiation} with
respect to $t$. The $2$ functions $x_{n}\left( t\right) $ are the dependent
variables, and other dependent variables $y_{n}\left( t\right) $ are
introduced below. Often below the dependence of these variables on $t$ shall
\textit{not} be \textit{explicitly} displayed, when this omission is
unlikely to cause any misunderstanding. The $6$ ($t$-independent) parameters
$c_{n\ell }$ are \textit{a priori} \textit{arbitrary}, but \textit{a
posteriori} we shall identify $2$ \textit{constraints} on their values; and
other $t$-independent parameters---such as $a_{nm},$ $b_{nm}$, etc.---shall
be introduced below. All variables and parameters can be \textit{complex}
numbers (but of course the subcase in which they are \textit{real}
numbers is of special interest in applicative contexts); we shall instead
generally think of the independent variable $t$ as \textit{time}, but
analytic continuation to \textit{complex} values of $t$---and of other
analogous time-like variables such as $\tau ,$ see below---shall also be
discussed. Generally each of the $2$ indices $n$ and $m$ take the $2$ values
$1$ and $2$, and the index $\ell $ the $3$ values $1,$ $2,$ $3$. $\ \
\blacksquare $

The system (\ref{1}) is a prototypical system of \textit{nonlinearly-coupled}
ODEs and as such has over time been studied in many theoretical
investigations and also utilized in an enormous number of applicative
contexts; a much too large research universe to make it possible to mention
all relevant references. Here we limit ourselves to quote the path-breaking
papers by Ren\'{e} Garnier \cite{RG5960}, and the very recent papers \cite%
{CCL2020} and \cite{CF2021}, whose topics are quite close to those treated
in the present paper, as discussed in the last two \textbf{Sections 6 }and
\textbf{7}, where possible future developments are also tersely outlined%
; and just one textbook reference \cite{DLA2006} (of
course the interested reader can trace additional references from those
quoted in these sources).

The main finding of the present paper is the identification (see \textbf{%
Sections 2}, \textbf{3} and\textbf{\ 4}) of a subclass of the model (\ref{1}%
)---characterized by $2$\textit{\ explicit algebraic constraints} on the $6$
coefficients $c_{n\ell }$ (see below the $2$ eqs. (\ref{Cons}))---which then
allows the \textit{explicit} solution of the \textit{initial-values problem}
for this system (\ref{1}), as detailed in \textbf{Proposition 2-2}.

Invariance properties of the system (\ref{1}) and some simplifications of it
are reported in \textbf{Section 5}.

Some extensions of the model (\ref{1}) to analogous systems with \textit{%
non-homogeneous quadratic} right-hand sides---including \textit{isochronous}
versions--- are discussed in \textbf{Section 6}.

A comparison with previous findings, and a very terse mention of possible
future developments, are provided in \textbf{Section 7}.

Let us complete this introductory \textbf{Section 1} with a terse
review---complementing the analogous treatment provided in \cite{CCL2020}%
---of the various meanings that can be given to the property of a system of
nonlinear ODEs to be \textit{solvable}, and more specifically to be \textit{%
explicitly solvable}.

As already noted in \cite{CCL2020}, the statement that a system of nonlinear
ODEs---such as (\ref{1})---is \textit{solvable by quadratures} is somewhat
\textit{misleading}, when it only implies that the independent variable $t$
can be identified as a function of an appropriate combination of the
dependent variables represented by an \textit{integral} which cannot be
\textit{explicitly} performed or that can be expressed as a \textit{named}
function---such as, say, a hypergeometric function---which cannot be readily
inverted. A less unsatisfactory outcome is when that function is a \textit{%
polynomial}, implying that its inversion yields an \textit{algebraic}
function, since this has significant implications, especially in terms of
the \textit{analytic} structure of the solution when considered as a
function of \textit{complex} $t$; although of course a \textit{generic}
polynomial cannot be \textit{explicitly} inverted---i. e., its \textit{roots}
identified---unless its degree does \textit{not} exceed $4$.

In the present paper the statement that a system of nonlinearly-coupled ODEs
is \textit{explicitly solvable} indicates that the solution of the
corresponding \textit{initial-values} problem can be exhibited as an \textit{%
elementary function} of the independent variable $t$, involving parameters
themselves expressed, in terms of the original parameters of the model, by
\textit{explicit} formulas only involving \textit{elementary} functions; the
final formulas expressing the parameters of the solution being nevertheless,
possibly, \textit{quite complicated}, being produced by a finite (generally
short) chain of \textit{explicit} relations applied sequentially (see
examples below).

\bigskip

\section{Main results}

The following $2$ Propositions are proven in the following \textbf{Section 3}%
.

\textbf{Proposition 2-1}. The \textit{explicit solution} of the \textit{%
initial-values} problem for the system
\begin{subequations}
\label{y12dot}
\begin{equation}
\dot{y}_{1}\left( t\right) =\left[ y_{1}\left( t\right) \right] ^{2}~,
\label{y1dot}
\end{equation}%
\begin{equation}
\dot{y}_{2}\left( t\right) =\rho _{1}\left[ y_{1}\left( t\right) \right]
^{2}+\rho _{2}y_{1}\left( t\right) y_{2}\left( t\right) +\left[ y_{2}\left(
t\right) \right] ^{2}~,  \label{y2dot}
\end{equation}%
where $\rho _{1}$ and $\rho _{2}$ are $2$ \textit{arbitrary} parameters,
reads as follows:
\end{subequations}
\begin{subequations}
\label{y12t}
\begin{equation}
y_{1}\left( t\right) =y_{1}\left( 0\right) \left[ 1-y_{1}\left( 0\right) t%
\right] ^{-1}~,  \label{y1t}
\end{equation}%
\begin{equation}
y_{2}\left( t\right) =y_{1}\left( 0\right) \left[ 1-y_{1}\left( 0\right) t%
\right] ^{-1}~u\left( t\right) ~,  \label{y2t}
\end{equation}%
\begin{equation}
u\left( t\right) =\frac{u_{+}\left[ u\left( 0\right) -u_{-}\right] -u_{-}%
\left[ u\left( 0\right) -u_{+}\right] \left[ 1-y_{1}\left( 0\right) t\right]
^{-\Delta }}{u\left( 0\right) -u_{-}-\left[ u\left( 0\right) -u_{+}\right] %
\left[ 1-y_{1}\left( 0\right) t\right] ^{-\Delta }}~,  \label{ut}
\end{equation}%
\begin{equation}
u\left( 0\right) =y_{2}\left( 0\right) /y_{1}\left( 0\right) ~,~~~u_{\pm
}=\left( 1-\rho _{2}\pm \Delta \right) /2~,  \label{u0pm}
\end{equation}%
\begin{equation}
\Delta =\sqrt{\left( 1-\rho _{2}\right) ^{2}-4\rho _{1}}~.  \label{Delta}
\end{equation}

This solution is valid for \textit{arbitrary initial} data $y_{1}\left(
0\right) $ and $y_{2}\left( 0\right) ,$ provided $y_{1}\left( 0\right) \neq
0 $. If instead $y_{1}\left( 0\right) =0$ implying $y_{1}\left( t\right)
=y_{1}\left( 0\right) =0$---in which case some of the formulas (\ref{y12t})
become undetermined---then of course $y_{2}\left( t\right) =y_{2}\left(
0\right) \left[ 1-y_{2}\left( 0\right) t\right] ^{-1}$ (see (\ref{y2dot})
with $y_{1}\left( t\right) =0$). $\ \ \blacksquare $

\textbf{Remark 2-1}. Note that this solution is clearly invariant under the
assignment of the sign of $\Delta $ (not defined by eq. (\ref{Delta})): see (%
\ref{ut}) and the definition (\ref{u0pm}) of the $2$ parameters $u_{\pm }$. $%
\ \ \blacksquare $

\textbf{Proposition 2-2}. The \textit{initial-values} problem---with
\textit{generic} initial data---for the system (\ref{1}) is \textit{%
explicitly solvable }provided the $6$ \textit{a priori arbitrary} parameters
$c_{n\ell }$ ($n=1,2;$ $\ell =1,2,3$) are expressed in terms of the $6=2+4$
\textit{a priori arbitrary} parameters $\rho _{1}$, $\rho _{2}$ and $a_{nm}$
\textit{or} $b_{nm}$ ($n=1,2;$ $m=1,2$) by the following formulas:
\end{subequations}
\begin{subequations}
\label{2cnm}
\begin{equation}
c_{n1}=b_{n1}\left( a_{11}\right) ^{2}+b_{n2}\left[ \rho _{1}\left(
a_{11}\right) ^{2}+\left( \rho _{2}a_{11}+a_{21}\right) a_{21}\right]
~,~~n=1,2~,
\end{equation}%
\begin{equation}
c_{n2}=2b_{n1}a_{11}a_{12}+b_{n2}\left[ 2\rho _{1}a_{11}a_{12}+\rho
_{2}\left( a_{11}a_{22}+a_{12}a_{21}\right) +2a_{21}a_{22}\right] ~,~~n=1,2,
\end{equation}%
\begin{equation}
c_{n3}=b_{n1}\left( a_{12}\right) ^{2}+b_{n2}\left[ \rho _{1}\left(
a_{12}\right) ^{2}+\left( \rho _{2}a_{12}+a_{22}\right) a_{22}\right]
~,~~n=1,2~.
\end{equation}%
Here the $4$ parameters $a_{nm}$ and the $4$ parameters $b_{nm}$ are related
by the following $4$ formulas:
\end{subequations}
\begin{subequations}
\label{abAB}
\begin{equation}
a_{11}=b_{22}/B~,~~a_{12}=-b_{12}/B~,~~a_{21}=-b_{21}/B~,~~a_{22}=b_{11}/B~,
\label{anm}
\end{equation}%
or, equivalently,%
\begin{equation}
b_{11}=a_{22}/A~,~~b_{12}=-a_{12}/A~,~~b_{21}=-a_{21}/A~,~~b_{22}=a_{11}/A~,
\label{bnm}
\end{equation}%
where%
\begin{equation}
A=a_{11}a_{22}-a_{12}a_{21}=B^{-1}~,  \label{A}
\end{equation}%
\begin{equation}
B=b_{11}b_{22}-b_{12}b_{21}=A^{-1}~;  \label{B}
\end{equation}%
obviously implying the possibility to express---via the formulas (\ref{2cnm}%
) with $n=1,2$---the $6$ coefficients $c_{n\ell }$ in terms of the $2$
parameters $\rho _{1}$, $\rho _{2}$ and \textit{either} the $4$ \textit{a
priori arbitrary} parameters $a_{nm}$ \textit{or} the $4$ \textit{arbitrary}
parameters $b_{nm}$.

Then the solution of the \textit{initial-values} problem for the system (\ref%
{1}) is related to the \textit{explicit} solution (\ref{y12t}) of the
corresponding \textit{initial-values} problem for the system (\ref{y12dot})
(see \textbf{Proposition 2-1}) via the following \textit{linear} relations:
\end{subequations}
\begin{subequations}
\label{y12x12}
\begin{equation}
y_{1}\left( t\right) =a_{11}x_{1}\left( t\right) +a_{12}x_{2}\left( t\right)
~,~~~y_{2}\left( t\right) =a_{21}x_{1}\left( t\right) +a_{22}x_{2}\left(
t\right) ~,  \label{y1y2x12}
\end{equation}%
\begin{equation}
x_{1}\left( t\right) =b_{11}y_{1}\left( t\right) +b_{12}y_{2}\left( t\right)
~,~~~x_{2}\left( t\right) =b_{21}y_{1}\left( t\right) +b_{22}y_{2}\left(
t\right) ~,  \label{x1x2y12}
\end{equation}%
which are easily seen to imply the relations (\ref{abAB}). $\ \ \blacksquare
$

\bigskip

\section{Proofs}

\subsection{Proof of \textbf{Proposition 2-1}}

In this subsection we provide for completeness a proof of \textbf{%
Proposition 2-1}, although this finding is rather elementary and by no means
new (see for instance \cite{RG5960}).

The fact that (\ref{y1t}) provides the solution of the \textit{initial-value}
problem for the ODE (\ref{y1dot}) is plain.

Hereafter we assume $y_{1}\left( 0\right) \neq 0$.

Then set
\end{subequations}
\begin{subequations}
\begin{equation}
y_{2}\left( t\right) =u\left( t\right) y_{1}\left( t\right) ~,~~~u\left(
t\right) =y_{2}\left( t\right) /y_{1}\left( t\right) ~,
\end{equation}%
hence%
\begin{equation}
\dot{u}\left( t\right) =\left[ \dot{y}_{2}\left( t\right) -u\left( t\right)
\dot{y}_{1}\left( t\right) \right] /y_{1}\left( t\right) ~,
\end{equation}%
hence, via (\ref{y12dot}) and (\ref{y1t}),
\end{subequations}
\begin{equation}
\dot{u}\left( t\right) =y_{1}\left( 0\right) \left[ 1-y_{1}\left( 0\right) t%
\right] ^{-1}\left\{ \left[ u\left( t\right) \right] ^{2}+\left( \rho
_{2}-1\right) u\left( t\right) +\rho _{1}\right\} ~,
\end{equation}%
and since clearly, via the definition of $u_{\pm }$ (see (\ref{u0pm})),%
\begin{equation}
\left[ u\left( t\right) \right] ^{2}+\left( \rho _{2}-1\right) u\left(
t\right) +\rho _{1}=\left[ u\left( t\right) -u_{+}\right] \left[ u\left(
t\right) -u_{-}\right] ~,
\end{equation}%
the ODE satisfied by $u\left( t\right) $ reads as follows:%
\begin{equation}
\dot{u}\left( t\right) \left\{ \left[ u\left( t\right) -u_{+}\right] \left[
u\left( t\right) -u_{-}\right] \right\} ^{-1}=y_{1}\left( 0\right) \left[
1-y_{1}\left( 0\right) t\right] ^{-1}~;
\end{equation}%
hence, again via the definition of $u_{\pm }$ (see (\ref{u0pm})),%
\begin{equation}
\dot{u}\left( t\right) \left\{ \left[ u\left( t\right) -u_{+}\right] ^{-1}-%
\left[ u\left( t\right) -u_{-}\right] ^{-1}\right\} =\Delta y_{1}\left(
0\right) \left[ 1-y_{1}\left( 0\right) t\right] ^{-1}  \label{21udot}
\end{equation}%
which can be immediately integrated, yielding (\ref{ut}).

The expression of $y_{2}\left( t\right) $ (see (\ref{y12t})) is thereby
validated, completing thereby the proof of \textbf{Proposition 2-1}.

\textbf{Remark 2.1-1}. Of course the expression (\ref{ut}) of $u\left(
t\right) $ is valid for \textit{generic} values of the relevant parameters.
For the special \textit{initial} values $u\left( 0\right) =u_{\pm }$ it
yields the trivial result $u\left( t\right) =u\left( 0\right) $. For the
\textit{special} values of the parameters $\rho _{1}$ and $\rho _{2}$ such
that $\rho _{1}=\left( 1-\rho _{2}\right) ^{2}/4$ implying $\Delta =0$ (see (%
\ref{Delta})) and (see (\ref{u0pm}))%
\begin{equation}
u_{+}=u_{-}=\bar{u}=\left( 1-\rho _{2}\right) /2~,  \label{upEqum}
\end{equation}%
the expression (\ref{ut}) of $u\left( t\right) $ is replaced by the
following formula (implied by (\ref{21udot}) with (\ref{upEqum})):%
\begin{equation}
u\left( t\right) =\frac{u\left( 0\right) +\bar{u}\left[ u\left( 0\right) -%
\bar{u}\right] ~\ln \left[ 1-y_{1}\left( 0\right) t\right] }{1+\left[
u\left( 0\right) -\bar{u}\right] ~\ln \left[ 1-y_{1}\left( 0\right) t\right]
}~.~~~\blacksquare
\end{equation}

\bigskip

\subsection{Proof of Proposition 2-2}

Let us $t$-differentiate the relations (\ref{x1x2y12}) respectively (\ref%
{y1y2x12}), getting
\begin{subequations}
\begin{equation}
\dot{x}_{1}=b_{11}\dot{y}_{1}+b_{12}\dot{y}_{2}~,~~~\dot{x}_{2}=b_{21}\dot{y}%
_{1}+b_{22}\dot{y}_{2}~,  \label{xyndot}
\end{equation}%
respectively%
\begin{equation}
\dot{y}_{1}=a_{11}\dot{x}_{1}+a_{12}\dot{x}_{2}~,~~~\dot{y}_{2}=a_{21}\dot{x}%
_{1}+a_{22}\dot{x}_{2}~.  \label{yxndot}
\end{equation}%
Hence, from the first of these $2$ pairs of relations, we get, via (\ref%
{y12dot}),
\end{subequations}
\begin{equation}
\dot{x}_{n}=b_{n1}\left( y_{1}\right) ^{2}+b_{n2}\left[ \rho _{1}\left(
y_{1}\right) ^{2}+\rho _{2}y_{1}y_{2}+\left( y_{2}\right) ^{2}\right]
~,~~~n=1,2~;
\end{equation}%
and then, via (\ref{y1y2x12}) and a bit of trivial algebra, the system (\ref%
{1}) with the expressions (\ref{2cnm}) of the coefficients $c_{n\ell }$.
\textbf{Proposition 2-2} is thereby proven.

\bigskip

\section{Inversion of the equations (\protect\ref{2cnm}) with (\protect\ref%
{abAB})}

In this Section we discuss the important problem to invert the system of
\textit{algebraic} equations (\ref{2cnm}) with (\ref{abAB}), i. e. to
express the $2$ parameters $\rho _{1}$, $\rho _{2}$ and the $4$ parameters $%
a_{nm}$---or, equivalently (see (\ref{abAB})), the $4$ parameters $b_{nm}$%
---in terms of the $6$ parameters $c_{n\ell }$; and we find $2$ \textit{%
constraints} on the $6$ parameters $c_{n\ell }$ which are required in order
to fulfill this task, hence are \textit{necessary} for the \textit{explicit}
solvability of the system (\ref{1}) via \textbf{Proposition 2.2}.

As a first step, let us note that the system of $2$ ODEs (\ref{yxndot})
implies, via the system (\ref{1}), the following $2$ ODEs:
\begin{subequations}
\begin{eqnarray}
\dot{y}_{n}=\left( a_{n1}c_{11}+a_{n2}c_{21}\right) \left( x_{1}\right)
^{2}+\left( a_{n1}c_{12}+a_{n2}c_{22}\right) x_{1}x_{2} &&  \notag \\
+\left( a_{n1}c_{13}+a_{n2}c_{23}\right) \left( x_{2}\right)
^{2}~,~~~n=1,2~, &&
\end{eqnarray}%
hence, via (\ref{x1x2y12}), the following system of $2$ ODEs:%
\begin{eqnarray}
\dot{y}_{n}=\left( a_{n1}c_{11}+a_{n2}c_{21}\right) \left(
b_{11}y_{1}+b_{12}y_{2}\right) ^{2} &&  \notag \\
+\left( a_{n1}c_{12}+a_{n2}c_{22}\right) \left(
b_{11}y_{1}+b_{12}y_{2}\right) \left( b_{21}y_{1}+b_{22}y_{2}\right) &&
\notag \\
+\left( a_{n1}c_{13}+a_{n2}c_{23}\right) \left(
b_{21}y_{1}+b_{22}y_{2}\right) ^{2}~,~~~n=1,2~, &&
\end{eqnarray}%
hence
\end{subequations}
\begin{equation}
\dot{y}_{n}=\gamma _{n1}\left( y_{1}\right) ^{2}+\gamma
_{n2}y_{1}y_{2}+\gamma _{n3}\left( y_{2}\right) ^{2}~,~~~n=1,2~,
\end{equation}%
with
\begin{subequations}
\begin{eqnarray}
\gamma _{n1}=\left( a_{n1}c_{11}+a_{n2}c_{21}\right) \left( b_{11}\right)
^{2}+\left( a_{n1}c_{12}+a_{n2}c_{22}\right) b_{11}b_{21} &&  \notag \\
+\left( a_{n1}c_{13}+a_{n2}c_{23}\right) \left( b_{21}\right)
^{2}~,~~~n=1,2~, &&
\end{eqnarray}%
\begin{eqnarray}
\gamma _{n2}=2\left( a_{n1}c_{11}+a_{n2}c_{21}\right) b_{11}b_{12}+\left(
a_{n1}c_{12}+a_{n2}c_{22}\right) \left( b_{11}b_{22}+b_{12}b_{21}\right) &&
\notag \\
+2\left( a_{n1}c_{13}+a_{n2}c_{23}\right) b_{21}b_{22}~,~~~n=1,2~, &&
\end{eqnarray}%
\begin{eqnarray}
\gamma _{n3}=\left( a_{n1}c_{11}+a_{n2}c_{21}\right) \left( b_{12}\right)
^{2}+\left( a_{n1}c_{12}+a_{n2}c_{22}\right) b_{12}b_{22} &&  \notag \\
+\left( a_{n1}c_{13}+a_{n2}c_{23}\right) \left( b_{22}\right)
^{2}~,~~~n=1,2~. &&
\end{eqnarray}%
And now a comparison of this system of ODEs with the system (\ref{y12dot})
implies
\end{subequations}
\begin{equation}
\gamma _{11}=\gamma _{23}=1~,~~~\gamma _{12}=\gamma _{13}=0~;~~\gamma
_{2n}=\rho _{n},~~n=1,2~,
\end{equation}%
hence
\begin{subequations}
\begin{eqnarray}
\gamma _{11}=\left( a_{11}c_{11}+a_{12}c_{21}\right) \left( b_{11}\right)
^{2}+\left( a_{11}c_{12}+a_{12}c_{22}\right) b_{11}b_{21} &&  \notag \\
+\left( a_{11}c_{13}+a_{12}c_{23}\right) \left( b_{21}\right) ^{2}=1~, &&
\label{gamma11}
\end{eqnarray}%
\begin{eqnarray}
\gamma _{12}=2\left( a_{11}c_{11}+a_{12}c_{21}\right) b_{11}b_{12}+\left(
a_{11}c_{12}+a_{12}c_{22}\right) \left( b_{11}b_{22}+b_{12}b_{21}\right) &&
\notag \\
+2\left( a_{11}c_{13}+a_{12}c_{23}\right) b_{21}b_{22}=0~, &&
\label{gamma12}
\end{eqnarray}

\begin{eqnarray}
\gamma _{13}=\left( a_{11}c_{11}+a_{12}c_{21}\right) \left( b_{12}\right)
^{2}+\left( a_{11}c_{12}+a_{12}c_{22}\right) b_{12}b_{22} &&  \notag \\
+\left( a_{11}c_{13}+a_{12}c_{23}\right) \left( b_{22}\right) ^{2}=0~, &&
\label{gamma13}
\end{eqnarray}

\begin{eqnarray}
\gamma _{23}=\left( a_{21}c_{11}+a_{22}c_{21}\right) \left( b_{12}\right)
^{2}+\left( a_{21}c_{12}+a_{22}c_{22}\right) b_{12}b_{22} &&  \notag \\
+\left( a_{21}c_{13}+a_{22}c_{23}\right) \left( b_{22}\right) ^{2}=1~, &&
\label{gamma23}
\end{eqnarray}%
\begin{eqnarray}
\gamma _{21}=\left( a_{21}c_{11}+a_{22}c_{21}\right) \left( b_{11}\right)
^{2}+\left( a_{21}c_{12}+a_{22}c_{22}\right) b_{11}b_{21} &&  \notag \\
+\left( a_{21}c_{13}+a_{22}c_{23}\right) \left( b_{21}\right) ^{2}=\rho
_{1}~, &&  \label{gamma21}
\end{eqnarray}%
\begin{eqnarray}
\gamma _{22}=2\left( a_{21}c_{11}+a_{22}c_{21}\right) b_{11}b_{12}+\left(
a_{21}c_{12}+a_{22}c_{22}\right) \left( b_{11}b_{22}+b_{12}b_{21}\right) &&
\notag \\
+2\left( a_{21}c_{13}+a_{22}c_{23}\right) b_{21}b_{22}=\rho _{2}~; &&
\label{gamma22}
\end{eqnarray}%
namely, by setting,
\end{subequations}
\begin{equation}
\alpha _{n\ell }=a_{n1}c_{1\ell }+a_{n2}c_{2\ell }~,~~~n=1,2,~~\ell =1,2,3~,
\label{alphacnel}
\end{equation}%
the following $6$ equations:
\begin{subequations}
\label{alphanm}
\begin{equation}
\alpha _{11}\left( b_{11}\right) ^{2}+\alpha _{12}b_{11}b_{21}+\alpha
_{13}\left( b_{21}\right) ^{2}=1~,
\end{equation}%
\begin{equation}
2\alpha _{11}b_{11}b_{12}+\alpha _{12}\left(
b_{11}b_{22}+b_{12}b_{21}\right) +2\alpha _{13}b_{21}b_{22}=0~,
\end{equation}

\begin{equation}
\alpha _{11}\left( b_{12}\right) ^{2}+\alpha _{12}b_{12}b_{22}+\alpha
_{13}\left( b_{22}\right) ^{2}=0~,
\end{equation}

\begin{equation}
\alpha _{21}\left( b_{12}\right) ^{2}+\alpha _{22}b_{12}b_{22}+\alpha
_{23}\left( b_{22}\right) ^{2}=1~,
\end{equation}%
\begin{equation}
\alpha _{21}\left( b_{11}\right) ^{2}+\alpha _{22}b_{11}b_{21}+\alpha
_{23}\left( b_{21}\right) ^{2}=\rho _{1}~,
\end{equation}%
\begin{equation}
2\alpha _{21}b_{11}b_{12}+\alpha _{22}\left(
b_{11}b_{22}+b_{12}b_{21}\right) +2\alpha _{23}b_{21}b_{22}=\rho _{2}~.
\end{equation}

\textbf{Remark 4-1}. From the first $3$ of these relations---summing the
first and the third and summing or subtracting the second---we get the
following $2$ relations
\end{subequations}
\begin{equation}
\alpha _{11}\left( b_{11}\pm b_{12}\right) ^{2}+\alpha _{12}\left( b_{11}\pm
b_{12}\right) \left( b_{21}\pm b_{22}\right) +\alpha _{13}\left( b_{21}\pm
b_{22}\right) ^{2}=1~,
\end{equation}%
and summing and subtracting these $2$ relations we get the $2$ relations
\begin{subequations}
\begin{equation}
\alpha _{11}\left[ \left( b_{11}\right) ^{2}+\left( b_{12}\right) ^{2}\right]
+\alpha _{12}\left( b_{11}b_{21}+b_{12}b_{22}\right) +\alpha _{13}\left[
\left( b_{21}\right) ^{2}+\left( b_{22}\right) ^{2}\right] =1~,
\end{equation}

\begin{equation}
2\alpha _{11}b_{11}b_{12}+\alpha _{12}\left(
b_{11}b_{22}+b_{12}b_{21}\right) +2\alpha _{13}b_{21}b_{22}=0~.
\end{equation}

But we shall not use these formulas below. $\ \ \blacksquare $

Solving the first $3$ of the $6$ eqs. (\ref{alphanm}) we get the following
formulas for the $3$ quantities $\alpha _{1\ell }$ ($\ell =1,2,3$):
\end{subequations}
\begin{equation}
\alpha _{11}=\left( b_{22}\right) ^{2}/B^{2}~,~~\alpha
_{12}=-2b_{12}b_{22}/B^{2}~,~~\alpha _{13}=\left( b_{12}\right) ^{2}/B^{2}~,
\label{alpha1el}
\end{equation}%
and likewise solving the last $3$ of the $6$ eqs. (\ref{alphanm}) we get the
following formulas for the $3$ quantities $\alpha _{2\ell }$ ($\ell =1,2,3$%
):
\begin{subequations}
\label{alpha2el}
\begin{equation}
\alpha _{21}=\left[ \left( b_{21}\right) ^{2}+\left( b_{22}\right) ^{2}\rho
_{1}-b_{21}b_{22}\rho _{2}\right] /B^{2}~,
\end{equation}%
\begin{equation}
\alpha _{22}=-\left[ 2b_{11}b_{21}+2b_{12}b_{22}\rho _{1}-\left(
b_{11}b_{22}+b_{12}b_{21}\right) \rho _{2}\right] /B^{2}~,
\end{equation}%
\begin{equation}
\alpha _{23}=\left[ \left( b_{11}\right) ^{2}+\left( b_{12}\right) ^{2}\rho
_{1}-b_{11}b_{12}\rho _{2}\right] /B^{2}~;
\end{equation}%
of course, above and below, $B$ is defined in terms of the $4$ parameters $%
b_{nm}$ by eq. (\ref{B}).

Next, using the definitions (\ref{alphacnel}) and the relations (\ref{anm}),
we get the following $6$ algebraic equations, which only involve the $6$
parameters $\rho _{1},~\rho _{2}$ and $b_{nm}$ as well as the $6$ parameters
$c_{n\ell }$:
\end{subequations}
\begin{subequations}
\label{bnmrho12}
\begin{equation}
\left( b_{22}c_{11}-b_{12}c_{21}\right) B=\left( b_{22}\right) ^{2}~,
\label{bnma}
\end{equation}%
\begin{equation}
\left( b_{22}c_{12}-b_{12}c_{22}\right) B=-2b_{12}b_{22}~,  \label{bnmb}
\end{equation}%
\begin{equation}
\left( b_{22}c_{13}-b_{12}c_{23}\right) B=\left( b_{12}\right) ^{2}~,
\label{bnmc}
\end{equation}%
\begin{equation}
\left( -b_{21}c_{11}+b_{11}c_{21}\right) B=\left( b_{21}\right) ^{2}+\left(
b_{22}\right) ^{2}\rho _{1}-b_{21}b_{22}\rho _{2}~,  \label{bnmrho12d}
\end{equation}%
\begin{equation}
\left( b_{21}c_{12}-b_{11}c_{22}\right) B=2b_{11}b_{21}+2b_{12}b_{22}\rho
_{1}-\left( b_{11}b_{22}+b_{12}b_{21}\right) \rho _{2}~,  \label{bnmrho12e}
\end{equation}%
\begin{equation}
\left( -b_{21}c_{13}+b_{11}c_{23}\right) B=\left( b_{11}\right) ^{2}+\left(
b_{12}\right) ^{2}\rho _{1}-b_{11}b_{12}\rho _{2}~.  \label{bnmrho12f}
\end{equation}%
In all these formulas $B$ is of course again defined in terms of the $4$
parameters $b_{nm}$ by the formula (\ref{B}).

Solving for $\rho _{1}$ and $\rho _{2}$ the $2$ linear eqs. (\ref{bnmrho12d}%
) and (\ref{bnmrho12e}) we get
\end{subequations}
\begin{subequations}
\begin{eqnarray}
\rho _{1} &=&\left\{ \left( b_{21}\right) ^{2}\left(
1-b_{12}c_{11}-b_{22}c_{12}\right) +\left( b_{11}\right)
^{2}b_{22}c_{21}\right.  \notag \\
&&\left. +b_{11}b_{21}\left[ b_{12}c_{21}-b_{22}\left( c_{11}-c_{22}\right) %
\right] \right\} /\left( b_{22}\right) ^{2}~,  \label{rho1a}
\end{eqnarray}%
\begin{equation}
\rho _{2}=\left(
2b_{21}-2b_{12}b_{21}c_{11}-b_{21}b_{22}c_{12}+2b_{11}b_{12}c_{21}+b_{11}b_{22}c_{22}\right) /b_{22}~;
\label{rho2a}
\end{equation}%
likewise, solving the $2$ linear eqs. (\ref{bnmrho12e}) and (\ref{bnmrho12f}%
), we get
\end{subequations}
\begin{subequations}
\label{rho}
\begin{eqnarray}
\rho _{1} &=&\left\{ b_{12}\left( b_{21}\right) ^{2}c_{13}+b_{11}b_{21}\left[
b_{22}c_{13}+b_{12}\left( c_{12}-c_{23}\right) \right] \right.  \notag \\
&&\left. +\left( b_{11}\right) ^{2}\left( 1-b_{12}c_{22}-b_{22}c_{23}\right)
\right\} /\left( b_{12}\right) ^{2}~,  \label{rho1b}
\end{eqnarray}%
\begin{equation}
\rho _{2}=\left[ b_{11}\left( 2-b_{12}c_{22}-2b_{22}c_{23}\right)
+b_{12}b_{21}c_{12}+2b_{21}b_{22}c_{13}\right] /b_{12}~;  \label{rho2b}
\end{equation}%
and likewise, solving the $2$ linear eqs. (\ref{bnmrho12f}) and (\ref%
{bnmrho12d}), we get
\end{subequations}
\begin{subequations}
\begin{equation}
\rho _{1}=\left[ \left( b_{21}\right) ^{2}b_{22}c_{13}+\left( b_{11}\right)
^{2}b_{12}c_{21}+b_{11}b_{21}\left( 1-b_{12}c_{11}-b_{22}c_{23}\right) %
\right] /(b_{12}b_{22})~,  \label{rho1c}
\end{equation}%
\begin{equation}
\rho _{2}=b_{21}/b_{22}+\left[ b_{12}(-b_{21}c_{11}+b_{11}c_{21})\right]
/b_{22}+\left( b_{11}+b_{21}b_{22}c_{13}-b_{11}b_{22}c_{23}\right) /b_{12}~.
\label{rho2c}
\end{equation}

Any one of these $3$ pairs of formulas provides an \textit{explicit}
expression of the $2$ parameters $\rho _{1}$ and $\rho _{2}$ in terms of the
$4$ parameters $b_{nm}$ and the $6$ parameters $c_{n\ell }.$ Hence hereafter
we may only focus on the problem to express the $4$ parameters $b_{nm}$ in
terms of the $6$ parameters $c_{n\ell }$.

Indeed, by identifying $2$ different expressions of the parameter $\rho _{1}$
or $\rho _{2}$ as given just above, we obtain additional formulas involving
\textit{only} the $4$ parameters $b_{nm}$ and the $6$ parameters $c_{n\ell
}. $ In particular by identifying the $2$ expressions (\ref{rho2a}) and (\ref%
{rho2b}) we get the following formula:
\end{subequations}
\begin{subequations}
\label{4bb12}
\begin{eqnarray}
&&b_{12}\left(
2b_{21}-2b_{12}b_{21}c_{11}-b_{21}b_{22}c_{12}+2b_{11}b_{12}c_{21}+b_{11}b_{22}c_{22}\right)
\notag \\
&=&b_{22}\left[ b_{11}\left( 2-b_{12}c_{22}-2b_{22}c_{23}\right)
+b_{12}b_{21}c_{12}+2b_{21}b_{22}c_{13}\right] ~;  \label{bb1}
\end{eqnarray}%
and likewise by identifying the $2$ expressions (\ref{rho1b}) and (\ref%
{rho1c}) we get the following formula:%
\begin{eqnarray}
&&b_{22}\left\{ b_{11}b_{21}\left[ b_{22}c_{13}+b_{12}\left(
c_{12}-c_{23}\right) \right] +\left( b_{11}\right) ^{2}\left(
1-b_{12}c_{22}-b_{22}c_{23}\right) \right\}  \notag \\
&=&b_{12}\left[ \left( b_{11}\right) ^{2}b_{12}c_{21}+b_{11}b_{21}\left(
1-b_{12}c_{11}-b_{22}c_{23}\right) \right] ~.  \label{bb2}
\end{eqnarray}

Our final task is to extract as much information as possible on the
dependence of the $4$ parameters $b_{nm}$ on the $6$ parameters $c_{n\ell }$%
, from these $2$ equations (\ref{4bb12}) and from the $3$ eqs. (\ref{bnma}),
(\ref{bnmb}), (\ref{bnmc}), or rather from $2$ of their $3$ ratios, which
clearly read as follows:
\end{subequations}
\begin{subequations}
\label{ratios}
\begin{equation}
-2b_{12}\left( b_{22}c_{11}-b_{12}c_{21}\right) =b_{22}\left(
b_{22}c_{12}-b_{12}c_{22}\right) ~,  \label{rat1}
\end{equation}%
\begin{equation}
b_{12}\left( b_{22}c_{12}-b_{12}c_{22}\right) =-2b_{22}\left(
b_{22}c_{13}-b_{12}c_{23}\right) ~,  \label{rat2}
\end{equation}%
\begin{equation}
\left( b_{12}\right) ^{2}\left( b_{22}c_{11}-b_{12}c_{21}\right) =\left(
b_{22}\right) ^{2}\left( b_{22}c_{13}-b_{12}c_{23}\right) ~;  \label{rat3}
\end{equation}%
each one of these $3$ formulas (\ref{ratios}) is of course implied by the
other $2$.

Let us now introduce the auxiliary variable
\end{subequations}
\begin{equation}
\beta =b_{12}/b_{22}~.  \label{Defbeta}
\end{equation}%
Then, by dividing the $2$ eqs.(\ref{rat1}) and (\ref{rat2}) by $\left(
b_{22}\right) ^{2}$ we get the following $2$ \textit{quadratic} equations
for this quantity:
\begin{subequations}
\label{Eqbeta}
\begin{equation}
2c_{21}\beta ^{2}+\left( c_{22}-2c_{11}\right) \beta -c_{12}=0~,
\label{Eqbetaa}
\end{equation}%
\begin{equation}
c_{22}\beta ^{2}-\left( c_{12}-2c_{23}\right) \beta -2c_{13}=0~.
\label{Eqbetab}
\end{equation}%
Subtracting the \textit{second} of these $2$ eqs. multiplied by $2c_{21}$
from the \textit{first} itself multiplied by $c_{22}$ we get a first-degree
equation for $\beta $, the solution of which reads
\end{subequations}
\begin{equation}
\beta =\left( c_{12}c_{22}-4c_{13}c_{21}\right) /\left[ \left(
c_{22}-2c_{11}\right) c_{22}+2c_{21}\left( c_{12}-2c_{23}\right) \right] ~;
\label{beta}
\end{equation}%
and inserting this determination of $\beta $ in the $2$ eqs. (\ref{Eqbeta})
we finally get the following $2$ \textit{explicit constraints} on the $6$
coefficients $c_{n\ell }$:
\begin{subequations}
\label{Cons}
\begin{eqnarray}
2c_{21}\left( c_{12}c_{22}-4c_{13}c_{21}\right) ^{2} &&  \notag \\
+\left( c_{22}-2c_{11}\right) \left( c_{12}c_{22}-4c_{13}c_{21}\right) \left[
\left( c_{22}-2c_{11}\right) c_{22}+2c_{21}\left( c_{12}-2c_{23}\right) %
\right]  &&  \notag \\
-c_{12}\left[ \left( c_{22}-2c_{11}\right) c_{22}+2c_{21}\left(
c_{12}-2c_{23}\right) \right] ^{2}=0~, &&  \label{Cona}
\end{eqnarray}%
\begin{eqnarray}
c_{22}\left( c_{12}c_{22}-4c_{13}c_{21}\right) ^{2} &&  \notag \\
-\left( c_{12}-2c_{23}\right) \left( c_{12}c_{22}-4c_{13}c_{21}\right) \left[
\left( c_{22}-2c_{11}\right) c_{22}+2c_{21}\left( c_{12}-2c_{23}\right) %
\right]  &&  \notag \\
-2c_{13}\left[ \left( c_{22}-2c_{11}\right) c_{22}+2c_{21}\left(
c_{12}-2c_{23}\right) \right] ^{2}=0~. &&  \label{Conb}
\end{eqnarray}%
These $2$ \textit{constraints} on the $6$ coefficients $c_{n\ell }$ must be
satisfied in order that the \textit{initial-values problem} of the system (%
\ref{1}) be \textit{explicitly solvable} as detailed by \textbf{Proposition
2-2}. Note that each of these \textit{constraints} is a \textit{quintic}
algebraic equation for the $6$ coefficients $c_{n\ell }$; but eq. (\ref{Cona}%
) is only \textit{quadratic} for $c_{11},$ $c_{13},$ $c_{23}$ and \textit{%
cubic} for $c_{12},$ $c_{21}$, $c_{22}$; while eq. (\ref{Conb}) is only
\textit{quadratic} for $c_{11},$ $c_{13},$ $c_{21},$ $c_{23},$ \textit{cubic}
for $c_{12}$ and \textit{quartic} for $c_{22}$.

\textbf{Remark 4.2.} The last sentence above suggests the most
convenient approaches to be employed in order to evaluate the implications
of the $2$ \textit{constraints}\ (\ref{Cons}) in the
special cases---generally relevant in applicative contexts---when the $6$
coefficients $c_{n\ell }$ are all \textit{real} numbers. \
$\blacksquare $

Let us now complete the task of this Section, to express the $4$ parameters $%
b_{nm}$---hence as well the $4$ parameters $a_{nm}$: see (\ref{anm}) with (%
\ref{B})---in terms of the $6$ coefficients $c_{n\ell }$. Since the
definition (\ref{Defbeta}) of $\beta $ clearly implies
\end{subequations}
\begin{equation}
b_{12}=\beta b_{22}~,  \label{bbeta}
\end{equation}%
inserting this relation in the $3$ eqs. (\ref{bnma}), (\ref{bb1}) and (\ref%
{bb2}), we get the following $3$ algebraic equations:
\begin{subequations}
\label{bbb}
\begin{equation}
\left( c_{11}-\beta c_{21}\right) \left( b_{11}-\beta b_{21}\right) =1~,
\label{bbb1}
\end{equation}%
\begin{eqnarray}
\left[ 1-b_{22}\left( \beta ^{2}c_{21}+\beta c_{22}+c_{23}\right) \right]
b_{11} &&  \notag \\
+\left[ -\beta +b_{22}\left( \beta ^{2}c_{11}+\beta c_{12}+c_{13}\right) %
\right] b_{21} &=&0~,  \label{bbb2}
\end{eqnarray}%
\begin{eqnarray}
b_{21}b_{22}\left[ c_{13}+\beta \left( c_{12}-c_{23}\right) \right] +\left(
b_{11}\right) ^{2}\left[ 1-b_{22}\left( \beta c_{22}+c_{23}\right) \right] &&
\notag \\
-\beta \left\{ \beta b_{11}b_{22}c_{21}+b_{21}\left[ 1-b_{22}\left( \beta
c_{11}+c_{23}\right) \right] \right\} =0~. &&  \label{bbb3}
\end{eqnarray}

The first $2$ of these $3$ eqs. (\ref{bbb}) are a \textit{linear} system for
the $2$ quantities $b_{11}$ and $b_{22}$, which can be immediately solved
yielding
\end{subequations}
\begin{subequations}
\label{Solb1121}
\begin{equation}
b_{11}=B_{110}+\beta b_{21}~,  \label{Solb11}
\end{equation}%
\begin{equation}
b_{22}=1/\left( B_{220}+B_{221}b_{21}\right) ~,  \label{Solb22}
\end{equation}%
with%
\begin{equation}
B_{110}=1/\left( c_{11}-\beta c_{21}\right) ~,  \label{b110}
\end{equation}%
\begin{equation}
B_{220}=c_{23}+\beta c_{22}+\beta ^{2}c_{21}~,  \label{b220}
\end{equation}%
\begin{equation}
B_{221}=-\left( c_{11}-\beta c_{21}\right) \left[ c_{13}+\beta \left(
c_{12}-c_{23}\right) +\beta ^{2}\left( c_{11}-c_{22}\right) -c_{21}\beta ^{3}%
\right] ~.  \label{b221}
\end{equation}%
Note that, since $\beta $ is \textit{explicitly} expressed in terms of the $%
6 $ coefficients $c_{n\ell }$ (see (\ref{beta})), the $2$ formulas (\ref%
{Solb1121})---together with (\ref{bbeta})---provide \textit{explicit}
expressions of the $3$ parameters $b_{11},$ $b_{22},$ $b_{12}$ in terms of
the $6$ coefficients $c_{n\ell }$ and the parameter $b_{21}$. Thus to
complete our task we must express also this parameter $b_{21}$ in terms of
the $6$ coefficients $c_{n\ell }$. This task can be fulfilled by solving the
\textit{algebraic} equation (\ref{bbb3}), which---after the replacement of
the $3$ parameters $b_{11},$ $b_{22},$ $b_{21}$ via their expressions (\ref%
{bbeta}) and (\ref{Solb1121}) in terms of $b_{21}$ and the $6$ parameters $%
c_{n\ell }$---only features the still unknown parameter $b_{21}$ (of course
in addition to the $6$ parameters $c_{n\ell }$). And it can be shown---via
elementary if tedious calculations, which can be checked by \textbf{%
Mathematica}---that the eq. (\ref{bbb3}) takes then the following form:
\end{subequations}
\begin{subequations}
\label{CCC}
\begin{equation}
C_{0}+C_{1}b_{21}+C_{2}\left( b_{21}\right) ^{2}+C_{3}\left( b_{21}\right)
^{3}=0~,  \label{Eqb21}
\end{equation}%
with%
\begin{equation}
C_{0}=\beta ^{2}c_{21}\left( 1-c_{11}+\beta c_{21}\right) ~,  \label{CC0}
\end{equation}%
\begin{eqnarray}
C_{1}=\left( c_{11}-\beta c_{21}\right) \left\{ \beta ^{3}c_{21}-\left(
1-c_{11}+\beta c_{21}\right) \cdot \right. &&  \notag \\
\left. \cdot \left[ c_{13}+\beta \left( c_{12}-c_{23}\right) +\beta
^{2}\left( c_{11}-c_{22}\right) -2\beta ^{3}c_{21}\right] \right\} ~, &&
\label{CC1}
\end{eqnarray}

\begin{eqnarray}
&&C_{2}=-\beta \left( c_{11}-\beta c_{21}\right) ^{2}\left\{ c_{13}+\beta
\left( c_{12}-c_{23}\right) +\beta ^{2}\left( c_{11}-c_{22}\right) -2\beta
^{3}c_{21}\right.  \notag \\
&&\left. +\left( 1-c_{11}+\beta c_{21}\right) \left[ c_{13}-\beta \left(
c_{23}-c_{12}\right) +\beta ^{2}\left( c_{11}-c_{22}\right) -\beta ^{3}c_{21}%
\right] \right\} ~,  \label{CC2}
\end{eqnarray}%
\begin{equation}
C_{3}=-\beta ^{2}\left( c_{11}-\beta c_{21}\right) ^{3}\left[ c_{13}+\beta
\left( c_{12}-c_{23}\right) +\beta ^{2}\left( c_{11}-c_{22}\right) -\beta
^{3}c_{21}\right] ~.  \label{CC3}
\end{equation}

Since these $4$ coefficients $C_{k}$ ($k=0,1,2,3$) are all \textit{explicitly%
} expressed---via these formulas: see (\ref{beta}) and (\ref{CCC})---in
terms of the $6$ coefficients $c_{n\ell }$, it seems that to complete our
task all that still needs to be done is to solve the \textit{cubic} equation
(\ref{Eqb21}), which can of course be \textit{explicitly} solved via the
Cardano formulas.

But the situation is a bit more tricky, and in fact more simple.

The point is that, as we know, the $6$ parameters $c_{n\ell }$ cannot be
assigned freely; the success of the entire treatment requires that they
satisfy the $2$ \textit{constraints} (\ref{Cons}); and, as it happens, this
requirement seems to imply that the coefficient $C_{3}$ vanishes, $C_{3}=0$.
We have been unable to \textit{prove} this result \textit{explicitly}: note
that the expression of $C_{3}$ in terms of the $6$ coefficients $c_{n\ell }$
is quite complicated, also due to the complicated dependence of $\beta $ on
the coefficients $c_{n\ell }$ (see (\ref{beta})); and the $2$ \textit{%
constraints} (\ref{Cons}) are as well fairly complicated. But quite
convincing evidence of this fact is provided by the numerical examples
reported below, see \textbf{Subsection 4.1}.

Hence, the third-degree equation (\ref{Eqb21}) can be replaced by the
second-degree equation
\end{subequations}
\begin{subequations}
\begin{equation}
C_{0}+C_{1}b_{21}+C_{2}\left( b_{21}\right) ^{2}=0~\ ,  \label{4QuadrEqb12}
\end{equation}%
the $2$ solutions of which read of course as follows:%
\begin{equation}
b_{21}=b_{21\pm }\equiv \left( -C_{1}\pm \sqrt{\left( C_{1}\right)
^{2}-4C_{0}C_{2}}\right) /2~.  \label{4b12pm}
\end{equation}

This finding seems to complete our task to determine---in terms of $6$
coeffcients $c_{n\ell }$, \textit{arbitrarily} assigned except for the
requirement to satisfy the $2$ constraints (\ref{Cons})--- the $6$
parameters $\rho _{n}$ and $b_{nm}$ ($n,m=1,2$): see (\ref{beta}), (\ref{rho}%
), (\ref{bbeta}), (\ref{Solb1121}) and (\ref{4b12pm}). Hence to provide the
\textit{explicit} solution of the \textit{initial-values} problem of the
system (\ref{1}), as detailed by \textbf{Proposition 2-2} in terms of the $6$
parameters $\rho _{n}$ and $b_{nm}$.

But a doubt should still linger in the mind of the alert reader: the
solution of the \textit{initial-values} problem of the system (\ref{1})
should be \textit{unique}; but we just found $2$ \textit{different} values
for the parameter $b_{21}$ (see (\ref{4b12pm})), hence as well for the other
$3$ parameters $b_{nm}$ (see the eqs. (\ref{Solb1121}) and (\ref{bbeta})),
and as well for the parameters $\rho _{n}$ (see the eqs. (\ref{rho})). This
means that, if our treatment is correct, these different values must end up
yielding the \textit{same} solution for the variables $x_{n}\left( t\right) $%
. This "miracle" is indeed validated by a check of many specific examples,
as reported in the following \textbf{Subsection 4.1}; with the added
observation that---as implied by the observation that the eqs. (\ref{b221})
and (\ref{CC3}) clearly imply $B_{221}=C_{3}\left[ \beta \left( c_{11}-\beta
c_{21}\right) \right] ^{-2}$---we may conclude that the vanishing of the
parameter $C_{3}$ also implies the vanishing of the parameter $B_{221}$: $%
B_{221}=0$; implying (via (\ref{b220})) the replacement of the expression (%
\ref{Solb22}) of $b_{22}$ by the simple expression
\end{subequations}
\begin{subequations}
\label{4b22b12}
\begin{equation}
b_{22}=1/\left( c_{23}+\beta c_{22}+\beta ^{2}c_{21}\right) ~,  \label{4b22}
\end{equation}%
and as a consequence also the replacement of (\ref{bbeta}) with%
\begin{equation}
b_{12}=\beta /\left( c_{23}+\beta c_{22}+\beta ^{2}c_{21}\right) ~.
\label{4b12}
\end{equation}%
These simpler formulas expressing the $2$ parameters $b_{22}$ and $b_{12}$
directly via the parameters $c_{n\ell }$ (recall (\ref{beta})) imply that
the values of these $2$ parameters are not affected by the $2$-valued
indeterminacy affecting the other $2$ parameters $b_{11}$ and $b_{21}$ (see (%
\ref{Solb11}) and (\ref{4b12pm})) as well as the values of the $2$
parameters $\rho _{n}$ (see (\ref{rho})).

\bigskip

\subsection{Specific solvable examples}

Let us introduce this Subsection by emphasizing that---due to the \textit{%
explicit} character of the formulas (\ref{2cnm}) expressing the $6$
coefficients $c_{n\ell }$ in terms of the $6$ parameters $\rho _{n}$ and $%
b_{nm}$ (or, equivalently, $a_{nm}$: see (\ref{abAB}))---it is quite easy to
manufacture examples of the system (\ref{1}) which are \textit{explicitly
solvable} via our treatment: all one has to do is to input an \textit{%
arbitrary} assignment of these $6$ parameters $\rho _{n}$ and $b_{nm}$ in
these formulas (\ref{2cnm}).

In this Subsection we report only $3$ examples of the system (\ref{1}) which
are \textit{explicitly solvable} via the technique described in the present
paper. But we also tested several other such examples, which are not
reported here; they all confirmed the assertion (that $C_{3}=0$) mentioned
in the last part of \textbf{Section 4}. Of course it shall be likewise easy
for the interested reader to identify in this manner other systems (\ref{1})
\textit{explicitly solvable} via the technique introduced in this paper (see
\textbf{Propositions 2- 2 }and \textbf{2-1}).

By inserting the values of the parameters $c_{n\ell }$---obtained by the
simple procedure described in the first paragraph of this Subsection---in
the relevant formulas written above (see \textbf{Section 4}), we verified
that they of course do satisfy the $2$ constraints (\ref{Cons}); that they
always do yield a vanishing value for the parameter $C_{3}$ (and also for
the parameter $B_{221}$); we obtained specific values for each of the $2$
parameters $b_{22}$ and $b_{12}$ (of course, the same as those originally
employed to determine the set of coefficients $c_{n\ell }$); while we
obtained instead $2$ alternative determinations for the couple of parameters
$b_{11}$ and $b_{21}$ and also for the couple of parameters $\rho _{1}$ and $%
\rho _{2}$. And moreover---remarkably: although this "miracle" was
expected---we verified that these $2$ different determinations yield---via
the relevant formulas of \textbf{Proposition 2-2} and \textbf{2-1 }(see eqs.
(\ref{x1x2y12}), (\ref{y12t}), (\ref{4b22b12}), (\ref{4b12pm}), (\ref{Solb11}%
), (\ref{rho}))---the same, \textit{unique}, solution of the \textit{%
initial-values} problem of the system (\ref{1}).

The \textit{first} example is identified by the following assignments of the
$6$ coefficients $c_{n\ell }$:
\end{subequations}
\begin{subequations}
\label{41Example1}
\begin{equation}
c_{11}=7/3,~~c_{12}=2,~~c_{13}=3~,~~c_{21}=-1,~~c_{22}=-2~,~~~c_{23}=-3~.
\label{41cnel}
\end{equation}

The corresponding values of the parameters $\beta ,$ $b_{12}$ and $b_{22}$
are
\begin{equation}
\beta =-3~,~~~b_{12}=1/2~,~~~b_{22}=-1/6~,  \label{41beta}
\end{equation}%
while for the values of the parameters $b_{11},$ $b_{21},$ $\rho _{1},$ $%
\rho _{2}$ and $\Delta $ (see (\ref{Delta})) we get
\begin{equation}
b_{11}=0~,~~b_{21}=-1/2,~~\rho _{1}=3/2~,~~\rho _{2}=0~,~~\Delta =\mathbf{i}%
\sqrt{5}~,  \label{41rhoa}
\end{equation}%
or%
\begin{equation}
b_{11}=1~,~~b_{21}=-5/6,~~\rho _{1}=7/2~,~~\rho _{2}=4~,~~\Delta =\mathbf{i}%
\sqrt{5}~;  \label{41rhob}
\end{equation}%
note the equality of the $2$ determinations of the parameter $\Delta ,$
which are of course essential for the final outcome, namely the following
\textit{unique explicit} solution of the \textit{initial-values} problem of
the system (\ref{1}) with (\ref{41cnel}):
\begin{eqnarray}
x_{1}\left( t\right) =\mathbf{i}\left\{ 3\left[ x_{1}\left( 0\right)
+3x_{2}\left( 0\right) \right] \left[ (2-\mathbf{i}\sqrt{5})x_{1}\left(
0\right) +3x_{2}\left( 0\right) \right. \right. &&  \notag \\
\left. -\left. \left[ (2+\mathbf{i}\sqrt{5})x_{1}\left( 0\right)
+3x_{2}\left( 0\right) \right] \left[ 1+\left( 2/3\right) t(x_{1}\left(
0\right) +3x_{2}\left( 0\right) \right] ^{\mathbf{i}\sqrt{5}}\right]
\right\} /D_{1}\left( t\right) ~, &&
\end{eqnarray}

\begin{eqnarray}
x_{2}\left( t\right) =\left\{ 3\left[ x_{1}\left( 0\right) +3x_{2}\left(
0\right) \right] \left[ -3x_{1}\left( 0\right) -2x_{2}\left( 0\right) -%
\mathbf{i}\sqrt{5}x_{2}\left( 0\right) \right. \right. &&  \notag \\
\left. \left. +\left[ 3x_{1}\left( 0\right) +(2-\mathbf{i}\sqrt{5}%
)x_{2}\left( 0\right) \right] \left\{ 1+\left( 2/3\right) t\left[
x_{1}\left( 0\right) +3x_{2}\left( 0\right) \right] \right\} ^{\mathbf{i}%
\sqrt{5}}\right] \right\} /D_{1}\left( t\right) ~, &&  \notag \\
&&
\end{eqnarray}%
\begin{eqnarray}
D_{1}\left( t\right) =\left\{ 3+2t\left[ x_{1}\left( 0\right) +3x_{2}\left(
0\right) \right] \right\} \left[ (-7-\mathbf{i}\sqrt{5})x_{1}\left( 0\right)
+(-3-3\mathbf{i}\sqrt{5})x_{2}\left( 0\right) \right. &&  \notag \\
\left. +\left[ (7-\mathbf{i}\sqrt{5})x_{1}\left( 0\right) +3(1-\mathbf{i}%
\sqrt{5})x_{2}\left( 0\right) \right] \left\{ 1+\left( 2/3\right) t\left[
x_{1}\left( 0\right) +3x_{2}\left( 0\right) \right] \right\} ^{\mathbf{i}%
\sqrt{5}}\right] ~. &&  \notag \\
&&
\end{eqnarray}

The \textit{second} example is identified by the following assignments of
the $6$ coefficients $c_{n\ell }$:
\end{subequations}
\begin{subequations}
\label{41Example2}
\begin{equation}
c_{11}=c_{12}=c_{13}=1~,~~~c_{21}=1/8~,~~~c_{22}=2~,~~~c_{23}=-1~.
\label{42cnel}
\end{equation}

The corresponding data read then as follows:
\begin{equation}
\beta =2~,~~~b_{12}=4/7~,~~~b_{22}=2/7~,
\end{equation}%
and
\begin{equation}
b_{11}=0~,~~b_{21}=-2/3,~~\rho _{1}=7/9~,~~\rho _{2}=-4/3~,~~\Delta =\sqrt{%
7/3}~,
\end{equation}%
or%
\begin{equation}
b_{11}=1~,~~b_{21}=-1/6,~~\rho _{1}=-35/144~,~~\rho _{2}=13/6~,~~\Delta =%
\sqrt{7/3}~;
\end{equation}%
yielding the following \textit{unique explicit} solution of the \textit{%
initial-values} problem of the system (\ref{1}) with (\ref{42cnel}):%
\begin{eqnarray}
&x_{1}\left( t\right) =\left\{ 4\left[ x_{1}\left( 0\right) -2x_{2}\left(
0\right) \right] \left[ \left( 5-\sqrt{21}\right) x_{1}\left( 0\right)
+4x_{2}\left( 0\right) \right] \right. &  \notag \\
&\left. -\left[ \left( 5+\sqrt{21}\right) x_{1}\left( 0\right) +4x_{2}\left(
0\right) \right] \left\{ 1-\left( 3/4\right) \left[ x_{1}\left( 0\right)
-2x_{2}\left( 0\right) \right] t\right\} ^{\sqrt{7/3}}\right\} /D_{2}\left(
t\right) ~,&  \notag \\
&&
\end{eqnarray}

\begin{eqnarray}
&x_{2}\left( t\right) =\left\{ 4\left[ x_{1}\left( 0\right) -2x_{2}\left(
0\right) \right] \left[ -x_{1}\left( 0\right) -\left( 5+\sqrt{21}\right)
x_{2}\left( 0\right) \right] \right. &  \notag \\
&\left. +\left[ x_{1}\left( 0\right) +\left( 5-\sqrt{21}\right) x_{2}\left(
0\right) \right] \left\{ 1-\left( 3/4\right) \left[ x_{1}\left( 0\right)
-2x_{2}\left( 0\right) \right] t\right\} ^{\sqrt{7/3}}\right\} /D_{2}\left(
t\right) ~,&  \notag \\
&&
\end{eqnarray}%
\begin{eqnarray}
&&D_{2}\left( t\right) =\left\{ -4+3\left[ x_{1}\left( 0\right)
-2x_{2}\left( 0\right) \right] t\right\} \left\{ \left( -7+\sqrt{21}\right)
x_{1}\left( 0\right) -2\left( 7+\sqrt{21}\right) x_{2}\left( 0\right) \right.
\notag \\
&&\left. +\left[ \left( 7+\sqrt{21}\right) x_{1}\left( 0\right) +2\left( 7-%
\sqrt{21}\right) x_{2}\left( 0\right) \right] \left\{ 1-\left( 3/4\right) %
\left[ x_{1}\left( 0\right) -2x_{2}\left( 0\right) \right] t\right\} ^{\sqrt{%
7/3}}\right\} .  \notag \\
&&
\end{eqnarray}

The \textit{third} example is identified by the following assignments of the
$6$ coefficients $c_{n\ell }$:
\end{subequations}
\begin{subequations}
\label{41Example3}
\begin{eqnarray}
c_{11} &=&-19/169~,~~~c_{12}=-265/507~,~~~c_{13}=110/1521~,  \notag \\
c_{21} &=&-27/169~,~~~c_{22}=-1/169~,~~~c_{23}=-36/169~.
\end{eqnarray}

The corresponding data read then as follows:
\begin{equation}
\beta =5/3~,~~~b_{12}=-5/2~,~~~b_{22}=-3/2~,
\end{equation}%
and
\begin{equation}
b_{11}=0~,~~b_{21}=-39/10,~~\rho _{1}=-11/25~,~~\rho _{2}=17/10~,~~\Delta
=3/2~,
\end{equation}%
or%
\begin{equation}
b_{11}=1~,~~b_{21}=-33/10,~~\rho _{1}=-14/25~,~~\rho _{2}=9/10~,~~\Delta
=3/2~;
\end{equation}%
yielding the following \textit{unique explicit} solution of the \textit{%
initial-values} problem of the system (\ref{1}) with (\ref{42cnel}):%
\begin{eqnarray}
x_{1}\left( t\right) =\left\{ \left[ 3x_{1}\left( 0\right) -5x_{2}\left(
0\right) \right] \left\{ 3861x_{1}\left( 0\right) -858x_{2}\left( 0\right)
\right. \right. &&  \notag \\
\left. +78\left[ 9x_{1}\left( 0\right) +11x_{2}\left( 0\right) \right]
\left\{ 1-\left( 2/39\right) \left[ 3x_{1}\left( 0\right) -5x_{2}\left(
0\right) \right] t\right\} ^{3/2}\right\} /D_{3}\left( t\right) ~, &&
\end{eqnarray}

\begin{eqnarray}
x_{2}\left( t\right) =-\left\{ 9\left[ 3x_{1}\left( 0\right) -5x_{2}\left(
0\right) \right] \left\{ 351x_{1}\left( 0\right) -78x_{2}\left( 0\right)
\right. \right. &&  \notag \\
\left. -39\left[ 9x_{1}\left( 0\right) +11x_{2}\left( 0\right) \right]
\left\{ 1-\left( 2/39\right) \left[ 3x_{1}\left( 0\right) -5x_{2}\left(
0\right) \right] t\right\} ^{3/2}\right\} /D_{3}\left( t\right) ~, &&
\end{eqnarray}%
\begin{eqnarray}
D_{3}\left( t\right) =\left\{ 39-2\left[ 3x_{1}\left( 0\right) -5x_{2}\left(
0\right) \right] t\right\} \left\{ 702x_{1}\left( 0\right) -156x_{2}\left(
0\right) \right. &&  \notag \\
\left. -39\left[ 9x_{1}\left( 0\right) +11x_{2}\left( 0\right) \right]
\left\{ 1-\left( 2/39\right) \left[ 3x_{1}\left( 0\right) -5x_{2}\left(
0\right) \right] t\right\} ^{3/2}\right\} . &&
\end{eqnarray}

\bigskip

\section{Invariance property and simplifications}

In this short section we report for completeness a rather obvious invariance
property and some possible trivial simplifications of the system (\ref{1}).
They amount to the elementary observation that the $2$ dependent variables
\end{subequations}
\begin{subequations}
\label{xnhat}
\begin{equation}
\hat{x}_{n}\left( \tau \right) \equiv \left( \mu _{n}/\lambda \right)
~x_{n}\left( t\right) ~,~~~\hat{t}\equiv \lambda ~t~,  \label{landamu}
\end{equation}%
with $\lambda $ and $\mu _{n}$ \textit{a priori arbitrary nonvanishing}
parameters, satisfy---\textit{mutatis mutandis}---essentially the same
system (\ref{1}) as the $2$ dependent variables $x_{n}\left( t\right) $:
\begin{equation}
\hat{x}_{n}^{\prime }\left( \hat{t}\right) \equiv d~\hat{x}_{n}\left( \hat{t}%
\right) /d\hat{t}=\hat{c}_{n1}\left[ \hat{x}_{1}\left( \hat{t}\right) \right]
^{2}+\hat{c}_{n2}\hat{x}_{1}\left( \hat{t}\right) \hat{x}_{2}\left( \hat{t}%
\right) +\hat{c}_{n3}\left[ \hat{x}_{2}\left( \hat{t}\right) \right]
^{2}~,~~~n=1,2~,
\end{equation}%
with%
\begin{equation}
\hat{c}_{n1}=\mu _{n}\left( \mu _{1}\right) ^{-2}c_{n1}~,~~\hat{c}_{n2}=\mu
_{n}\left( \mu _{1}\mu _{2}\right) ^{-1}c_{n2}~,~~\hat{c}_{n3}=\mu
_{n}\left( \mu _{2}\right) ^{-2}c_{n3}~,~~n=1,2~.  \label{chat}
\end{equation}

For $\mu _{1}=\mu _{2}=1$ this property identifies the \textit{invariance }%
of the system (\ref{1}) under a simultaneous rescaling of the independent
and dependent variables: see (\ref{xnhat}).

\textbf{Remark 5-1}. Both \textit{constraints} (\ref{Cons}) are invariant
under the transformation (\ref{chat}). $\ \ \blacksquare $

The \textit{simplifications} correspond to the possibility to replace---by
an appropriate rescaling of dependent variables---$1$ of the $3$ parameters $%
c_{1\ell }$ and $1$ of the $3$ parameters $c_{2\ell }$ by an \textit{%
arbitrary} number (of course, nonvanishing; for instance, just \textit{unity}%
); thereby reducing the number of \textit{a priori arbitrary} coefficients $%
c_{n\ell }$ from $6$ to $4$. For instance the assignment
\end{subequations}
\begin{subequations}
\label{Semplif}
\begin{equation}
\mu _{1}=c_{11}~,~~~\mu _{2}=c_{23}~,  \label{Simpla}
\end{equation}%
implies%
\begin{eqnarray}
\hat{c}_{11} &=&1~,~~\hat{c}_{12}=c_{12}/c_{23}~,~~\hat{c}%
_{13}=c_{11}c_{13}\left( c_{23}\right) ^{-2}~,  \notag \\
\hat{c}_{21} &=&c_{21}c_{23}\left( c_{11}\right) ^{-2}~,~~\hat{c}%
_{22}=c_{22}/c_{11}~,~~\hat{c}_{23}=1~.  \label{Simplb}
\end{eqnarray}

\textbf{Remark 5-2}. The simplification (\ref{Semplif}) applied to the $3$
examples characterized by the assignments (\ref{41Example1}), (\ref%
{41Example2}) respectively (\ref{41Example3}) yields $3$ models (see (\ref%
{xnhat})) characterized by the following assignments of the coefficients $%
c_{n\ell }$:
\end{subequations}
\begin{equation}
\hat{c}_{12}=-2/3~,~~~\hat{c}_{13}=7/9~,~~~\hat{c}_{21}=27/49~,~~~\hat{c}%
_{22}=-6/7~,
\end{equation}%
respectively%
\begin{equation}
\hat{c}_{12}=-1~,~~~\hat{c}_{13}=1~,~~~\hat{c}_{21}=-1/8~,~~~\hat{c}_{22}=2~,
\end{equation}%
respectively%
\begin{equation}
\hat{c}_{12}=\frac{265}{108}~,~~~\hat{c}_{13}=-\frac{1045}{5832}~,~~~\hat{c}%
_{21}=\frac{972}{361}~,~~~\hat{c}_{22}=\frac{1}{19}~;
\end{equation}%
of course in all $3$ cases with $\hat{c}_{11}=\hat{c}_{23}=1$. $\ \
\blacksquare $

\bigskip

\section{Extensions and isochronous models}

In this Section we tersely outline some simple extensions of the system (\ref%
{1}) to the case with \textit{non-homogeneous quadratic} right-hand sides,
as well as some related systems obtained by a well-known change of
variables---see for instance \cite{C2008}---which allows the identification
of analogous systems featuring the remarkable property to be \textit{%
isochronous}.

An elementary way to extend the \textit{autonomous} system (\ref{1})
featuring ODEs with \textit{homogeneous quadratic} right-hand sides to an,
\textit{also autonomous}, system with \textit{non-homogeneous quadratic}
right-hand sides is via the following---easily \textit{invertible}---change
of \textit{independent} variables:
\begin{equation}
z_{n}\left( t\right) =\exp \left( \eta t\right) x_{n}\left( \tilde{t}\right)
+\bar{z}_{n}~,~~~\tilde{t}=\left[ \exp \left( \eta t\right) -1\right] /\eta
~,~~~n=1,2~,  \label{zn}
\end{equation}%
where the $3$ parameters $\bar{z}_{1}$, $\bar{z}_{2}$ and $\eta $ are
\textit{a priori arbitrary}. Thereby the system (\ref{1}) gets transformed
into the following system:
\begin{subequations}
\label{Eqzn}
\begin{eqnarray}
\dot{z}_{n}\left( t\right) =c_{n1}\left[ z_{1}\left( t\right) \right]
^{2}+c_{n2}z_{1}\left( t\right) z_{2}\left( t\right) +c_{n3}\left[
z_{2}\left( t\right) \right] ^{2} &&  \notag \\
+\eta z_{n}\left( t\right) +d_{n1}z_{1}\left( t\right) +d_{n2}z_{2}\left(
t\right) +d_{n3}~,~~~n=1,2~, &&  \label{zndot}
\end{eqnarray}%
with the $6$ "new" parameters $d_{n\ell }$ expressed in terms of the $6$
"old" parameters $c_{n\ell }$ and of the $3$ "new" parameters $\bar{z}_{1}$,
$\bar{z}_{2}$ and $\eta $ as follows:
\begin{eqnarray}
d_{n1} &=&-2c_{n1}\bar{z}_{1}-c_{n2}\bar{z}_{2}~,~~~d_{n2}=-2c_{n1}\bar{z}%
_{2}-c_{n2}\bar{z}_{1}~,  \notag \\
d_{n3} &=&-\eta \bar{z}_{n}+c_{n1}\left( \bar{z}_{1}\right) ^{2}+c_{n2}\bar{z%
}_{1}\bar{z}_{2}+c_{n3}\left( \bar{z}_{2}\right) ^{2}~,~~~n=1,2~.  \label{dn}
\end{eqnarray}

Of course the \textit{solvability} properties of the original system (\ref{1}%
) carry over to the system (\ref{Eqzn}).

Moreover---if the \textit{solvability} of the system (\ref{zndot}), via (\ref%
{zn}) and \textbf{Propositions 2-2} and \textbf{2-1---}features a parameter $%
\Delta $ (see (\ref{Delta})) which is a \textit{real rational }number ($%
\Delta =k_{1}/k_{2}$ with $k_{1}$ an \textit{arbitrary integer} and $k_{2}$
an \textit{arbitrary positive integer}), then clearly the system (\ref{zndot}%
), with
\end{subequations}
\begin{equation}
\eta =\mathbf{i}\omega ~,  \label{etadnel}
\end{equation}%
---where $\mathbf{i}$ is the imaginary unit, $\mathbf{i}^{2}=-1$, and $%
\omega $ is an \textit{arbitrary nonvanishing real} number---features the
remarkable property to be \textit{isochronous}: namely \textit{all} its
solutions $z_{n}\left( t\right) $ are \textit{periodic} with the same period
$T=2\pi k_{2}/\left\vert \omega \right\vert $,
\begin{equation}
z_{n}\left( t+T\right) =z_{n}\left( t\right) ~,~~~n=1,2~.  \label{Iso}
\end{equation}%
Readers wondering about the validity of this---rather obvious: see (\ref{zn}%
), (\ref{y12x12}), (\ref{ut}) and (\ref{etadnel})---conclusion are advised
to have a look, for instance, at the book \cite{C2008}.

\textbf{Remark 6-1}. Of course the presence of the \textit{imaginary}
parameter $\eta =\mathbf{i}\omega $ in the right-hand side of the system (%
\ref{zndot}) with (\ref{etadnel}) implies that its solutions are necessarily
%\textit{complex}, $z_{n}\left( t\right) \equiv \func{Re}\left[ z_{n}\left(
%t\right) \right] +\mathbf{i}\func{Im}\left[ z_{n}\left( t\right) \right] $;
\textit{complex}, $z_{n}\left( t\right) \equiv \mbox{Re}\left[ z_{n}\left(
t\right) \right] +\mathbf{i}\mbox{Im}\left[ z_{n}\left( t\right) \right] $;
entailing a corresponding doubling, from $2$ to $4$, of the number of
nonlinearly-coupled ODEs for the \textit{real} version of this system,
satisfied by the $4$ \textit{real} dependent variables $\mbox{Re}\left[
z_{n}\left( t\right) \right] $ and $\mbox{Im}\left[ z_{n}\left( t\right) %
\right] $, $n=1,2$; and clearly in this case it would be natural to also
consider the $6$ parameters $c_{n\ell }$ (as well of course as the $6$
parameters $\rho _{n}$ and $b_{nm}$ related to them) and the $2$ parameters $%
\bar{z}_{n}$ to be themselves \textit{complex} numbers. $\ \ \blacksquare $

\textbf{Remark 6-2}. The interested reader might wish to compute the
relevant formulas for the \textit{isochronous} case associated to the
\textit{third example} reported in \textbf{Subsection 4.1}. \ \ $%
\blacksquare $

\bigskip

\section{Comparison with previous findings and outlook}

The system (\ref{1}) treated in this paper is identical to the system
treated in the recent paper \cite{CCL2020}; it is therefore appropriate to
compare the approach and the findings reported in that paper with those
reported in the present paper.

The methodologies used in \cite{CCL2020} and in the present paper have much
in common, but there is a significant difference. In the present paper we
started from the simpler, \textit{explicitly solvable} model (\ref{y12dot})
and we then investigated in which cases the general system (\ref{1}) with $6$
\textit{a priori arbitrary} coefficients $c_{n\ell }$ can be reduced---via a
\textit{time-independent linear} transformation of the $2$ dependent
variables, see (\ref{y12x12})---to the simpler, \textit{explicitly solvable}
system (\ref{y12dot}). We found that this is indeed possible, but only if
the $6$ \textit{a priori arbitrary} coefficients $c_{n\ell }$ satisfy the $2$
\textit{constraints} (\ref{Cons}). This allowed us to conclude that the
special subclass of the systems (\ref{1}) identified by these $2$ \textit{%
constraints} is \textit{explicitly solvable} in terms of \textit{elementary}
functions, and to display the solution of their \textit{initial-values}
problem.

The methodology employed in \cite{CCL2020} took as point of departure the
general system (\ref{1}) with $6$ arbitrary coefficients $c_{n\ell }$, but
then immediately proceeded to reduce it to a \textit{canonical form}%
---featuring at most only $2$ coefficients---via a \textit{time-independent
linear} transformation of the $2$ dependent variables (such as (\ref{y12x12}%
)); it then focussed on the discussion of the \textit{solvability} (by
quadratures) of those reduced systems, and moreover on the identification of
a specific subclass of such systems the solutions of which are \textit{%
algebraic}, i. e. identified as roots of \textit{explicitly time-dependent}
\textit{polynomials}. The procedure of reduction to \textit{canonical form}
is a bit complicated, but it has been shown by Fran\c{c}ois Leyvraz that the
first example treated in Subsection III.B of \cite{CCL2020} (see eqs. (38-41
there) is essentially equivalent---up to notational changes---to the model
treated in the present paper. We also take this opportunity to mention a
trivial misprint in eq. (10b) of \cite{CCL2020}, which identifies the
Newtonian equation $\ddot{\zeta}=\zeta ^{k}$ as \textit{algebraically
solvable} if $k=-(2n+1)/\left( 2n-1\right) $ or $k=-(n+1)/n$ with $n$ a
\textit{positive integer}: the first of these $2$ equalities should instead
read $k=-(n+2)/n$ yielding $k=-3,$ $-2,$ $-5/3,-3/2,...$ (note that the
values of $k$ yielded by the definition $k=-(2n+1)/\left( 2n-1\right) $ with
$n$ an \textit{arbitrary positive} \textit{integer} coincide with those
yielded by the definition $k=-(m+2)/m$ \textit{only} if $m$ is an \textit{%
odd positive integer}).

Let us conclude by expressing the wishful hope that the type of approach
used in the present paper be also applicable to other systems of nonlinear
ODEs or PDEs---possibly also with \textit{discrete} rather than \textit{%
continuous} independent variables.

\bigskip

\section{Acknowledgements}

It is a pleasure to thank our colleagues Robert Conte, Fran\c{c}ois Leyvraz
and Andrea Giansanti for very useful discussions. We also like to
acknowledge with thanks $2$ grants, facilitating our collaboration---mainly
developed via e-mail exchanges---by making it possible for FP to visit twice
the Department of Physics of the University of Rome "La Sapienza": one
granted by that University, and one granted jointly by the Istituto
Nazionale di Alta Matematica (INdAM) of that University and by the
International Institute of Theoretical Physics (ICTP) in Trieste in the
framework of the ICTP-INdAM "Research in Pairs" Programme. Finally, we also\
like to thank Fernanda Lupinacci who, in these difficult times---with
extreme efficiency and kindness---facilitated all the arrangements necessary
for the presence of FP with her family in Rome.

\label{lastpage}

\end{document}